\documentclass[12pt]{article}
\usepackage{amsmath,amssymb,amsthm}
\def\a{\alpha}
\def\b{\beta}

\def\e{\epsilon}

\def\g{\gamma}

\def\z{\zeta}
\def\th{\theta}

\def\n{\nu}

\def\r{\rho}

\def\om{\omega}

\def\E{{\cal E}}
\def\Gp{G_{n,p}}
\def\whp{{\bf whp}}
\def\whps{{\bf whp }}
\def\C{{\cal C}}
\def\E{{\cal E}}

\newtheorem{lemma}{Lemma}
\newtheorem{theorem}{Theorem}

\newcommand{\proofstart}{{\bf Proof\hspace{2em}}}

\newcommand{\proofend}{\hspace*{\fill}\mbox{$\Box$}}

\newcommand{\expect}{\mbox{\bf E}}

\newcommand{\rdown}[1]{\lfloor #1 \rfloor}
\newcommand{\rdup}[1]{\lceil #1 \rceil}
\def\Prob{{\bf Pr}}
\def\ki{\rdown{k/i}}

\begin{document}
\baselineskip 20pt \lineskip 3pt
\title{Covering the edges of a random graph by cliques}
\author{Alan Frieze\thanks{Department of Mathematics, Carnegie Mellon
University. Supported in part by NSF grants
CCR9024935 and CCR9225008.} \and  Bruce Reed\thanks{Department of Mathematics, Carnegie Mellon
University and Equipe Combinatoire, CNRS, Universit\'e Pierre et Marie Curie, Paris}}
\maketitle

\section{Introduction}
The {\em clique cover number} $\th_1(G)$ of a graph $G$ is the 
minimum number of cliques required to cover the edges of graph $G$. 
In this paper we consider $\th_1(\Gp)$, for $p$ constant. (Recall that 
in the random graph $\Gp$, each of the ${n\choose 2}$ edges occurs 
independently with probability $p$). Bollob\'as, Erd\H{o}s, Spencer 
and West \cite{BESW} proved that \whps (i.e. with probability 1-o(1) as $n\rightarrow\infty$)
$${(1-o(1))n^2\over 4(\log_2n)^2}\leq \th_1(G_{n,.5})\leq {cn^2\ln \ln n\over (\ln n)^2}.$$
They implicitly conjecture that the $\ln \ln n$ factor in the upper 
bound is unnecessary and in this paper we prove
\begin{theorem}
\label{th1}
There exist constants $c_i=c_i(p)>0,i=1,2$ such that \whp
$${c_1n^2\over (\ln n)^2}\leq \th_1(G_{n,p})\leq {c_2n^2\over (\ln n)^2}.$$
\end{theorem}

{\bf Remark 1:} a simple use of a martingale tail inequality shows 
that $\th_1$ is close to its mean with very high probability.
\section{Proof of Theorem \protect{\ref{th1}}}
We write $a_n\approx b_n$ if $a_n/b_n\rightarrow 1$ as $n\rightarrow\infty$.

The lower bound is simple as the number of edges $m$ of $\Gp$ \whps satisfies
$$m\approx {np^2\over 2}$$
and the size of the largest clique $\om=\om(\Gp)$ \whps satisfies
$$\om\approx2\log_bn$$
where $b=1/p$. We may thus choose $c_1\approx (\ln b)^2p/2.$

The upper bound requires more work. Our method does not seem to yield 
the correct value for $c_2$ and so we will not work hard to keep $c_2$ 
small. Let $\a$ be some small constant and let
$$k=\rdown{\a\log_bn}.$$
We consider an algorithm for randomly selecting cliques to cover the edges 
of $G=\Gp$. It bears some relation to part of the algorithm described in 
Pippenger and Spencer \cite{PS}. At iteration $i$ we randomly select cliques 
of size $k_i=\rdown{k/i}$ none of whose edges are covered by previously chosen 
cliques. Our idea is to choose these cliques so that at the start of iteration 
$i$ the graph $G_i$ formed by the set $E_i$ of edges which have not been covered 
behaves, for our purposes, similarly to $G_{n,p_i}, p_i=pe^{1-i}$. That is it will 
contain about $m_i={n\choose 2}p_i$ edges, it will have about 
$N_i={n\choose k_i}p_i^{{k_i\choose 2}}$  cliques of size $k_i$ and 
the intersection of these cliques will be similar to that for the $k_i$-cliques 
in $G_{n,p_i}$. In particular, in both $G_{n,p_i}$ and $G_i$ almost all of the 
edges are in about $\z_i=N_i{k_i\choose 2}/m_i$ $k_i$-cliques.

Now in iteration $i$ we choose a set $\C_i$ of $k_i$-cliques from $G_i$ to add to 
our cover. The available cliques are chosen independently with probability about 
$1/\z_i$. By our assumptions on $G_i$, an edge is left uncovered with probability 
about $e^{-1}$. With a bit of care we can show that our assumptions continue to 
hold for $G_{i+1}$ as well.

We do this for $i_0=\rdup{4\ln\ln n}$ iterations. After this there are about 
${n\choose 2}pe(\ln n)^{-4}$ uncovered edges and we can add these as cliques 
of size two to the cover. In iteration $i$ we choose about 
$m_i/{k_i\choose 2}\approx n^2i^2pe^{1-i}k^{-2}$ cliques and so the total 
number of cliques used is $O(n^2/(\ln n)^2)$
as required.

We now need to describe our clique choosing process a little more formally: 
let $\C_{j,i}$ denote the set of $j$-cliques all of whose edges are in $E_i$. If
$$c_{s,j,i}={n-s\choose j-s}(be^i)^{{s\choose 2}-{j\choose 2}},$$
then $c_{s,j,i}$ is close to the expected number of cliques in $\C_{j,i}$ which 
contain a particular fixed clique in $\C_{s,i}$.

For a clique $S\in \C_{s,i}$ we let
$$X_{S,j,i}=|\{C\in \C_{j,i}:C\supseteq S\}|$$
and for integer $s\geq 0$,
$$X_{s,j,i}^*=\max\{X_{S,j,i}:S\in \C_{s,i}\}.$$
{\bf Algorithm} COVER

{\bf begin}

\hspace{2em} $E_1:=E(\Gp);\ \C_{COVER}:=\emptyset;$

\hspace{2em} {\bf for} $i$ =1 {\bf to} $i_0$ {\bf do}

\hspace{2em} begin

A: \hspace{3em} independently place each $C\in \C_{\rdown{k/i},i}$ into $\C_{COVER}$ with probability

\hspace{4em} ${X^*_{2,\ki,i}}^{-1}$;

B: \hspace{3em} for each $u\in E_i$ which is not covered by a clique in Step A, add $u$

\hspace{4em}(as a clique of size 2) to $\C_{COVER}$ with probability $\r_u$ where
$$e^{-1}-{X^*_2}^{-1}=\left(1-{1\over X^*_2}\right)^{X_u}(1-\r_u),$$
\hspace{4em}$X^*_2=X^*_{2,\rdown{k/i},i}$ and $X_u=X_{u,\rdown{k/i},i}$.

\hspace{2em} {\bf end}

\hspace{2em} $\C_{COVER}:=\C_{COVER}\cup E_{i_0+1}$.

{\bf end}

Observe first that the definition of $\r_u$ assumes that $X^*_2$ is large (which it is \whp) and so
\begin{eqnarray*}
\left(1-{1\over X^*_2}\right)^{X_u}&\geq&\left(1-{1\over X^*_2}\right)^{X^*_2}\\
&\geq & e^{-1}-{X^*_2}^{-1},
\end{eqnarray*}
and $\r_u$ is properly defined.

The following lemma contains the main core of the proof:
\begin{lemma}
\label{lem1}
Let $\E_i$ refer to the following two conditions:

(a)
$$X_{S,j,i}\leq (1+\b_i)c_{s,j,i},\hspace{1in}0\leq s\leq j\leq k/i\mbox{ and }S\in \C_{s,i},$$
where $\b_i=in^{-1/4}$,

(b)
$$X_{u,j,i}\geq (1-\g_i)c_{2,j,i},\hspace{1in}e\in E_i\mbox{ and }2\leq j\leq k/i$$
for all but at most $in^{31/16}$ edges, where $\g_i=in^{-1/16}$.

Then
\begin{eqnarray}
\label{eq1}
\Prob(\E_1)&=&1-o(n^{-1}),\\
\label{eq2}
\Prob(\E_{i+1}\mid \E_i) & \geq & 1-O(n^{-1/16}\log n).
\end{eqnarray}
\end{lemma}
We defer the proof of the lemma to the next section and show how to use it to 
prove Thereom \ref{th1}. Observe first that
\begin{equation}
\label{eq3}
{c_{s+1,j,i}\over c_{s,j,i}}=\left({j-s\over n-s}\right)(be^i)^s,
\end{equation}
and
\begin{equation}
\label{4}
c_{s,j,i}\geq n^{7/8}
\end{equation}
when $\a$ is small and $0\leq s< j\leq k/i.$

Next let $Y_i$ and $Z_i$ denote the number of $\rdown{k/i}$-cliques and edges 
respectively added to $\C_{COVER}$ in iteration $i$.
\begin{eqnarray}
\expect(Y_i\mid \E_i)&\leq&\expect
\left({X^*_{0,\ki,i}\over X^*_{2,\ki,i}}\bigg|\E_i\right)\nonumber\\
&\leq&(1+o(1)){c_{0,\ki,i}\over c_{2,\ki,i}}\nonumber\\
&\approx&{n^2i^2\over bk^2e^i},\label{eq4}
\end{eqnarray}
on using (\ref{eq3})

Since $Y_i$ is binomially distributed, we see using standard bounds 
on the tails of the binomial, that
$$\Prob\left(Y_i\geq {2n^2i^2\over bk^2e^i}\bigg|\E_i\right)\leq n^{-1}.$$
Thus
$$\Prob\left(\sum_{i=1}^{i_0}Y_i\geq \sum_{i=1}^{i_0}
{2n^2i^2\over bk^2e^i}\bigg|\E_0\right)=O\left({i_0\log n\over n^{1/16}}\right),$$
and so
\begin{equation}
\label{eq5}
\Prob\left(\sum_{i=1}^{i_0}Y_i\geq \sum_{i=1}^{i_0}{2n^2i^2\over bk^2e^i}\right)=o(1).
\end{equation}
Now a simple calculation gives
\begin{equation}
\label{eq6}
\r_u=O\left({X^*_2-X_u\over X^*_2}\right)
\end{equation}
and so
\begin{eqnarray*}
\expect(Z_i\mid \E_i)&=&O(in^{31/16}+\b_i|E_i|)\\
&=&O(n^{31/16}\ln n).
\end{eqnarray*}
Thus
$$\Prob(Z_i\geq n^{63/32}\mid \E_i)=O(n^{-1/32}\ln n)$$
and so
$$\Prob(\exists 1\leq i\leq i_0:Z_i\geq n^{63/32}\mid \E_0)=O(n^{-1/32}(\ln n)^2)$$
and
\begin{equation}
\label{eq7}
\Prob(\sum_{i=1}^{i_0}Z_i\geq i_0n^{63/32})=o(1).
\end{equation}
Also
\begin{eqnarray*}
\Prob(u\in E_{i+1}\mid u\in E_i)&=&\left(1-{1\over X^*_2}\right)^{X_u}(1-\r_u)\\
&<&e^{-1}.
\end{eqnarray*}
Thus
$$\expect(|E_{i_0+1}|)=O\left({n^2\over (\ln n)^4}\right)$$
and
\begin{equation}
\label{eq8}
\Prob\left(|E_{i_0+1}|\geq {n^2\over (\ln n)^3}\right)=o(1).
\end{equation}
Theorem \ref{th1} follows from (\ref{eq5}), (\ref{eq7}) and (\ref{eq8}) and
$$|\C_{COVER}|=\sum_{i=1}^{i_0}Y_i+\sum_{i=1}^{i_0}Z_i+|E_{i_0+1}|.$$
As we only use estimates for $X^*_{0,\ki,i}$ and $X^*_{2,\ki,i}$ the reader 
may wonder why it is necessary to prove Lemma \ref{lem1}(a) for $0\leq s\leq j\leq k/i$. 
The reason is simply that the lemma is proved by induction and we use a stronger induction 
hypothesis than the needed outcome.
\section{Proof of Lemma \protect{\ref{lem1}}}
If $s=j$ then $X_{S,j,i}=c_{s,j,i}=1$ and so we can assume $s<j$ from now on.

Let us first consider $\E_1$. Fix a set $S$ of size $s$, $0\leq s\leq k$. Assume it forms a 
clique in $G$. This does not condition any edges not contained in $S$. For a set $T$ let 
$N_c(T)$ denote the set of common neighbours of $T$ in $G$. We can enumerate the set of 
$j$-cliques containing $S$ as follows: choose $x_1\in N_c(S)$, 
$x_2\in N_c(S\cup \{x_1\}),\ldots,x_{j-s}\in N_c(S\cup \{x_1,x_2,\ldots,x_{j-s-1}\})$. 
The number of choices $\n_t$ for $x_t$ given $x_1,x_2,\ldots,x_{t-1}$ is distributed 
as $Bin(n-(s-t+1),p^{s+t-1})$. Thus for $0\leq \e\leq 1$
\begin{eqnarray*}
\Prob\left(\left|{\n_t\over (n-s-t+1)p^{s+t-1}}-1\right|\geq \e\right)&
\leq &2\exp\left\{-{\e^2(n-s-t+1)p^{s+t-1}\over 3}\right\}\\
&\leq &2\exp\{-\e^2n^{1-\a}/4\}.
\end{eqnarray*}
Putting $\e=n^{-1/3}$ we see that since there are $n^{O(\ln n)}$ choices for 
$x_1,x_2,\ldots ,x_{j-s}$,
$$\Prob\left(\left|{X_{S,j,0}\over c_{s,j,0}}-1\right|\geq n^{-1/3+o(1)}\right)
\leq \exp\{-n^{1/4}\}.$$
There are $n^{O(\ln n)}$ choices for $S$ and (\ref{eq1}) follows.

Assume now that $\E_i$ holds. We first prove
\begin{lemma}
\label{lem2}
Suppose $e_1,e_2,\ldots,e_t\in E_i$. Then
$$\Prob(e_t\in E_{i+1}\mid e_1,e_2,\ldots,e_{t-1}\in E_{i+1})=e^{-1}\left(1+
O\left({t\ln n\over n}\right)\right)$$
uniformly for $1\leq t\leq n^{1/2}$.
\end{lemma}
\proofstart
\begin{eqnarray}
\Prob(e_t\in E_{i+1}\mid e_1,e_2,\ldots,e_{t-1}\in E_{i+1})&\geq&\Prob(e_t\in E_{i+1})\label{n1}\\
&=&\left(1-{1\over X^*_2}\right)^{X_u}(1-\r_u)\nonumber\\
&=&e^{-1}-{X^*_2}^{-1}.\nonumber
\end{eqnarray}
Here $u=e_t,\ X_u=X_{u,\ki,i}$ and $X^*_2=X^*_{2,\ki,i}$ and inequality (\ref{n1}) 
follows from the fact that knowing $e_1,e_2,\ldots e_{t-1}\in E_{i+1}$ tells us that 
certain cliques (and edges) were not chosen for $\C_{COVER}$. On the other hand
\begin{eqnarray}
\Prob(e_t\in E_{i+1}\mid e_1,e_2,\ldots,e_{t-1}\in E_{i+1})&\leq&
\left(1-{1\over X^*_2}\right)^{X_u-tX^*_3}(1-\r_u)\label{11}\\
&=&(e^{-1}-{X^*_2}^{-1})\left(1-{1\over X^*_2}\right)^{tX^*_3}\nonumber\\
&=&e^{-1}\left(1+O\left({tX^*_3\over X^*_2}\right)\right),\nonumber
\end{eqnarray}
where $X^*_3=X^*_{3,\ki,i}$. If $\E_i$ holds then $X^*_3/X^*_2=O(\ln n/n)$.

Inequality (\ref{11}) follows from the fact that $e_t=u$ lies in at least $X_u-(t-1)X^*_3$ 
cliques which contain none of $e_1,e_2,\ldots,e_{t-1}$. This in turn arises from a two 
term inclusion-exclusion inequality and the fact that $e_t$ and $e_i$ together lie in at 
most $X_3^*$ cliques, for $1\leq i\leq t-1$.
\proofend

Now fix a set $S\in \C_{s,i}$ and let $X=X_{S,j,i+1}$ for some $j\leq k/(i+1)$. Condition 
on $S\in \C_{s,i+1}$. Let $\C_{S,j,i}=\{C\in \C_{j,i}:C\supseteq S\}$. Then
on using Lemma \ref{lem2}, we have
\begin{eqnarray}
\expect(X)&=&\sum_{C\in \C_{S,j,i}}\Prob(C\in \C_{j,i+1}\mid S\in \C_{s,i+1})\nonumber\\
&=&X_{S,j,i}\exp\left\{{s\choose 2}-{j\choose 2}\right\}\left(1+O\left({j^4\ln n
\over n}\right)\right),\label{z1}\\
&=&\expect(X_{S,j,0})\exp\left\{(i+1)\left({s\choose 2}-{j\choose 
2}\right)\right\}\left(1+O\left({j^4\ln n
\over n}\right)\right),\nonumber\\
\noalign{by induction on $i$}
&=&c_{s,j,0}\exp\left\{(i+1)\left({s\choose 2}-{j\choose
2}\right)\right\}\left(1+O\left({j^4\ln n
\over n}\right)\right),\nonumber\\
&=&c_{s,j,i+1}\left(1+O\left({j^4\ln n
\over n}\right)\right).\label{z1x}
\end{eqnarray}

We are going to use the Markov inequality
\begin{equation}
\label{z2}
\Prob(X\geq x)\leq {\expect((X)_r)\over (x)_r}
\end{equation}
where $(x)_r=x(x-1)(x-2)\ldots (x-r+1)$ and $r=\rdown{n^{3/8}}$.

Let ${\cal B}(\ell_2,\ell_3,\ldots,\ell_r)=\{(C_1,C_2,\ldots,C_r): (i)\ C_t\neq C_{t'}$ 
for $t\neq t'$, (ii) $C_t\in \C_{S,j,i}$, (iii) $|\C_t\cap(C_1\cup C_2\cup 
\cdots C_{t-1})|=s+\ell_t,$ for $t,t'=2,3,\ldots,r\}$. Then
$$\expect((X)_r)=\sum_{\ell_2,\ell_3,\ldots,\ell_r}\sum_{{\cal B}(\ell_2,\ell_3,\ldots,\ell_r)}
\Prob(C_1,C_2,\ldots,C_r\in \C_{j,i+1}\mid S\in \C_{s,i+1}).$$
From (\ref{z1})
$$\Prob(C_1\in \C_{j,i+1}|S\in \C_{s,i+1})=\exp\left\{{s\choose 2}-
{j\choose 2}\right\}\left(1+O\left({j^4\ln n\over n}\right)\right)$$
and
\begin{eqnarray*}
\Prob(C_t\in \C_{j,i+1}\mid C_1,C_2,\ldots,C_{t-1}\in \C_{j,i+1})&=
&\exp\left\{{s+\ell_t\choose 2}-{j\choose 2}\right\}\left(1+O\left({j^4\ln n\over n}\right)\right)\\
&=&\exp\left\{{s+\ell_t\choose 2}-{s\choose 2}\right\}{c_{s,j,i+1}
\over c_{s,j,i}}\left(1+O\left({j^4\ln n\over n}\right)\right)
\end{eqnarray*}
Also,
\begin{eqnarray*}
|{\cal B}(\ell_2,\ell_3,\ldots,\ell_r)|&\leq&\prod_{t=1}^r\left({(t-1)j-s
\choose \ell_t}X^*_{s+\ell_t,j,i}\right)\\
&\leq&\prod_{t=1}^r(rj)^{\ell_t}(1+\b_i)\left({b^{s+\ell_t}je^{i(s+\ell_t)}
\over n}\right)^{\ell_t}c_{s,j,i}.
\end{eqnarray*}
Hence,
\begin{eqnarray}
{\expect((X)_r)\over c_{s,j,i+1}^r}&\leq&\left(1+O\left({(\ln n)^4r
\over n}\right)\right)\sum_{\ell_2,\ell_3,\ldots,\ell_r}\prod_{t=1}^r(1+\b_i)
\left({e^{(\ell_t+2s-1)/2}rj^2(be^i)^{s+\ell_t}\over n}\right)^{\ell_t}\nonumber\\
&\leq&\left(1+O\left({(\ln n)^4r\over n}\right)\right)(1+\b_i)^r\sum_{\ell_2,\ell_3,\ldots,\ell_r}
\left({rk^2e^{3k}b^{2k}\over n}\right)^{\ell_2+\cdots+\ell_t}\\
&\leq&(1+rn^{-3/4})(1+\b_i)^r,\label{z3}
\end{eqnarray}
for $\a$ sufficiently small.

Hence, using (\ref{z2}),
\begin{eqnarray*}
\Prob(X\geq (1+\b_{i+1})c_{s,j,i+1})&\leq&{2(1+\b_i)^rc_{s,j,i+1}^r
\over ((1+\b_{i+1})c_{s,j,i+1})_r},\hspace{.5in}\mbox{by (\ref{z3})}\\
&\leq&3\left({1+\b_i\over 1+\b_{i+1}}\right)^r,\hspace{.75in}\mbox{using (\ref{4})}\\
&\leq&3\exp\left\{-{r(\b_{i+1}-\b_i)\over 1+\b_{i+1}}\right\}\\
&=&\exp\{-n^{1/8-o(1)}\}.
\end{eqnarray*}
There are $n^{O(\ln n)}$ choices for $S$ and $j$ and so part (a) of the lemma is 
proven. 

It remains only to deal with $X_{u,j,i+1}$ for an edge $u\in E_i$.
It follows from (\ref{z1x}) that if $X=X_{u,j,i+1}$ then
\begin{equation}
\label{z4}
\expect(X)=c_{2,j,i}\left(1+O\left({j^4\ln n
\over n}\right)\right),
\end{equation}
and from (\ref{z3}) that
\begin{equation}
\label{z5}
\expect(X(X-1))\leq\left(1+{3i\over n^{1/4}}\right)c_{2,j,i+1}^2.
\end{equation}
Suppose now that $X_{u,j,i}\geq (1-\g_i)c_{2,j,i}$. Then (\ref{z4}) and (\ref{z5}) imply that
\begin{align}
&\Prob(X\leq (1-\g_{i+1})c_{2,j,i+1})=\nonumber\\
&\Prob(\expect(X)-X\geq \expect(X)-(1-\g_{i+1})c_{2,j,i+1})\leq\nonumber\\
&\Prob\left(\expect(X)-X\geq (1-\g_i)c_{2,j,i}\exp\left\{1-{j\choose 2}\right\}
\left(1+O\left({j^4\ln n\over n}\right)\right)-(1-\g_{i+1})c_{2,j,i+1}\right)=\nonumber\\
&\Prob\left(\expect(X)-X\geq (1-\g_i)c_{2,j,i+1}
\left(1+O\left({j^4\ln n\over n}\right)\right)-(1-\g_{i+1})c_{2,j,i+1}\right)=\nonumber\\
&\Prob\left(\expect(X)-X\geq (1-o(1))n^{-1/16}c_{2,j,i+1}\right)\leq\nonumber\\
&\frac{(\expect(X)-X)^2}{(1-o(1))n^{-1/8}c^2_{2,j,i+1}}=\nonumber\\
&\frac{\expect(X(X-1))+\expect(X)-\expect(X)^2}{(1-o(1))n^{-1/8}c^2_{2,j,i+1}}\leq\nonumber\\
&\frac{\left(1+{3i\over n^{1/4}}\right)c_{2,j,i+1}^2+\left(1+O\left({j^4\ln n \over n}\right)\right)c_{2,j,i+1}
-c_{2,j,i+1}^2\left(1+O\left({j^4\ln n \over n}\right)\right)}
{(1-o(1))n^{-1/8}c^2_{2,j,i+1}}\le\nonumber\\
&\leq 4in^{-1/8}.\label{z6}
\end{align}
Now let $Z_{i+1}$ denote the number of edges $u\in E_{i+1}$ for which $X_{u,j,i+1}
\leq (1-\g_{i+1})c_{2,j,i+1}$ and $\hat{Z}_{i+1}$ those $u$ counted in $Z_{i+1}$ for 
which $X_{u,j,i}\geq (1-\g_i)c_{2,j,i}$. Then
$$Z_{i+1}\leq Z_i+\hat{Z}_{i+1}$$
and from (\ref{z6})
$$\expect(\hat{Z}_{i+1}\mid \E_i)\leq 4i|E_i|n^{-1/8}.$$
So
\begin{eqnarray*}
\Prob(Z_{i+1}\geq (i+1)n^{31/16}\mid \E_i)&\leq&\Prob(\hat{Z}_{i+1}\geq n^{31/16}\mid \E_i)\\
&=&O(n^{-1/16}\log n).
\end{eqnarray*}
this completes the proof of Lemma \ref{lem1}.

{\bf Acknowledgement:} We thank Reza Akhtar for pointing out an error in the paper.

\end{document}